\documentclass[12pt]{article}
\usepackage[english]{babel}
\usepackage{lineno} 
\usepackage{graphicx,multicol}
\usepackage{epic,eepic,epsfig}
\usepackage{caption}

\usepackage{amsmath,amsfonts,amsthm,amssymb}
\usepackage{pstricks,pstricks-add,pst-math,pst-xkey}

\usepackage{algorithm, algorithmic}

\usepackage{graphicx,pstricks,pst-node,amssymb}

\usepackage{tikz}

\pgfdeclarelayer{edgelayer}
\pgfdeclarelayer{nodelayer}
\pgfsetlayers{edgelayer,nodelayer,main}

\tikzstyle{none}=[inner sep=0pt]
\definecolor{hexcolor0xf81e1c}{rgb}{0.973,0.118,0.110}
\definecolor{hexcolor0x3c00ff}{rgb}{0.235,0.000,1.000}

\tikzstyle{whitevertex}=[circle,fill=white,draw=black, scale = 0.5]
\tikzstyle{redvertex}=[circle,fill=hexcolor0xf81e1c,draw=black, scale = 0.5]
\tikzstyle{bluevertex}=[circle,fill=hexcolor0x3c00ff,draw=black, scale = 0.5]
\tikzstyle{greenvertex}=[circle,fill=green,draw=black, scale=0.5]
\tikzstyle{purplevertex}=[circle,fill=magenta,draw=black, scale=0.5]
\tikzstyle{grayvertex}=[circle,fill=white,draw=gray, scale=0.5]
\tikzstyle{blackvertex}=[circle,fill=black,draw=black, scale=0.5]

\tikzstyle{textbox}=[rectangle,fill=none,draw=none]
\tikzstyle{box}=[rectangle,fill=none,draw=black]

\tikzstyle{arc}=[black, ->]
\tikzstyle{grayarc}=[gray, ->]
\tikzstyle{bluearc}=[blue, ->]
\tikzstyle{grayedge}=[draw=gray]
\tikzstyle{blueedge}=[draw=blue]
\tikzstyle{rededge}=[draw=red]
\tikzstyle{edge}=[draw=black]

\tikzstyle{vertex}=[circle, ,fill=white,draw=black, scale=0.5]

\tikzstyle{10circle}=[circle, scale=10.0,draw=black]
\tikzstyle{10oval}=[ellipse, scale=10.0,draw=black]

\usepackage{amsmath,amsfonts,amsthm,amssymb}
\usepackage{pstricks,pstricks-add,pst-math,pst-xkey}

\usepackage{algorithm, algorithmic}

\usepackage{graphicx,pstricks,pst-node,amssymb}

\usepackage{algorithm, algorithmic}

\begin{document}

\newtheorem{theorem}{\hspace{5mm}Theorem}[section]
\newtheorem{prp}[theorem]{\hspace{5mm}Proposition}
\newtheorem{definition}[theorem]{\hspace{5mm}Definition}
\newtheorem{lemma}[theorem]{\hspace{5mm}Lemma}
\newtheorem{corollary}[theorem]{\hspace{5mm}Corollary}
\newtheorem{alg}[theorem]{\hspace{5mm}Algorithm}
\newtheorem{sub}[theorem]{\hspace{5mm}Algorithm}
\newcommand{\induce}[2]{\mbox{$ #1 \langle #2 \rangle$}}
\newcommand{\2}{\vspace{2mm}}
\newcommand{\dom}{\mbox{$\rightarrow$}}
\newcommand{\ndom}{\mbox{$\not\rightarrow$}}
\newcommand{\compdom}{\mbox{$\Rightarrow$}}
\newcommand{\cdom}{\compdom}
\newcommand{\sdom}{\mbox{$\Rightarrow$}}
\newcommand{\lsd}{locally semicomplete digraph}
\newcommand{\lt}{local tournament}
\newcommand{\la}{\langle}
\newcommand{\ra}{\rangle}
\newcommand{\pf}{{\bf Proof: }}
\newtheorem{claim}{Claim}
\newcommand{\beq}{\begin{equation}}
\newcommand{\eeq}{\end{equation}}
\newcommand{\<}[1]{\left\langle{#1}\right\rangle}

\newcommand{\Z}{\mathbb{$Z$}}
\newcommand{\Q}{\mathbb{$Q$}}
\newcommand{\R}{\mathbb{$R$}}


\title{Chordality of locally semicomplete and weakly quasi-transitive digraphs}

\author{Jing Huang\ and\ Ying Ying Ye\thanks{Department of Mathematics and 
Statistics, University of Victoria, Victoria, B.C., Canada; huangj@uvic.ca 
(Research supported by NSERC)}}

\date{}

\maketitle

\begin{abstract}
Chordal graphs are important in the structural and algorithmic graph theory. 
A digraph analogue of chordal graphs was introduced by Haskin and Rose 
in 1973 but has not been a subject of active studies until recently when 
a characterization of semicomplete chordal digraphs in terms of forbidden 
subdigraphs was found by Meister and Telle.

Locally semicomplete digraphs, quasi-transitive digraphs, and extended semicomplete
digraphs are amongst the most popular generalizations of semicomplete digraphs. 
We extend the forbidden subdigraph characterization of semicomplete chordal 
digraphs to locally semicomplete chordal digraphs. We introduce a new class of 
digraphs, called weakly quasi-transitive digraphs, which contains quasi-transitive
digraphs, symmetric digraphs, and extended semicomplete digraphs, but is 
incomparable to the class of locally semicomplete digraphs. 
We show that weakly quasi-transitive digraphs can be recursively constructed by 
simple substitutions from transitive oriented graphs, semicomplete digraphs,
and symmetric digraphs. This recursive construction of weakly quasi-transitive 
digraphs, similar to the one for quasi-transitive digraphs, demonstrates
the naturalness of the new digraph class. As a by-product, we prove that the 
forbidden subdigraphs for semicomplete chordal digraphs are the same for weakly 
quasi-transitive chordal digraphs. The forbidden subdigraph
characterization of weakly quasi-transitive chordal digraphs generalizes not only
the recent results on quasi-transitive chordal digraphs and extended
semicomplete chordal digraphs but also the classical results on chordal graphs.
\end{abstract}

\section{Introduction}

We consider digraphs which do not contain loops or multiple arcs but may contain
digons (i.e., pairs of arcs joing vertices in opposite directions). If an arc
is contained in a digon then it is called a {\em symmetric arc}. A digraph which
does not contain a symmetric arc is called an {\em oriented graph}. A digraph which
contains only symmetric arcs is called a {\em symmetric digraph}. Graphs may be
viewed as symmetric digraphs. 

Two vertices in a digraph $D$ are {\em adjacent} and refered to as {\em neighbours}
of each other if there is at least one arc between them. We say that $u$ is an 
{\em in-neighbour} of $v$ or $v$ an {\em out-neighbour} of $u$ if $uv$ is an arc 
in $D$ (symmetric or not). The set of all in-neighbours of a vertex $v$ is 
denoted by $N^-(v)$ and the set of all out-neighbours of $v$ is denoted by $N^+(v)$.
We use $S(D)$ to denote the spanning subdigraph of $D$ whose arc set consists of 
all symmetric arcs in $D$.

A vertex $v$ in a digraph $D$ is {\em di-simplicial} if for every $u \in N^-(v)$
and $w \in N^+(v)$ with $u \neq w$, $uw$ is an arc of $D$. A digraph $D$ is 
{\em chordal} if every induced subdigraph of $D$ contains a di-simplicial vertex. 
It follows that every chordal digraph $D$ has a vertex ordering 
$v_1, v_2, \dots, v_n$ such that $v_i$ is a di-simplicial vertex in the subdigraph 
of $D$ induced by $v_i, v_{i+1}, \dots, v_n$ for each $i \geq 1$. Such an ordering 
is called a {\em perfect elimination ordering} of $D$. 

Perfect elimination orderings of digraphs arise in the study of sparse linear 
systems by Gaussian elimination, cf. \cite{hr}. When a digraph is symmetric, 
di-simplicial vertices coincide with simplicial vertices of its underlying graph. 
Thus, a symmetric digraph is chordal if and only if its underlying graph
is chordal. It is well-known that chordal graphs are precisely the graphs which do 
not contain an induced cycle of length $\geq 4$, cf. \cite{golumbic}. 

Little is known about the forbidden structure of chordal digraphs. In particular,
there is no known characterization of chordal digraphs by forbidden subdigraphs.
Recently, Meister and Telle \cite{meister} found a forbidden subdigraph 
characterization for semicomplete chordal digraphs. 
A digraph $D$ is {\em semicomplete} if between any two vertices there is at least
one arc. The following theorem is proved in \cite{meister}. 

\begin{theorem} \cite{meister} \label{basic} 
A semicomplete digraph $D$ is chordal if and only if $S(D)$ is chordal and $D$ does
not contain any of the digraphs in Figure \ref{Forbidden1} as an induced subdigraph.
\qed
\end{theorem}

\begin{center}
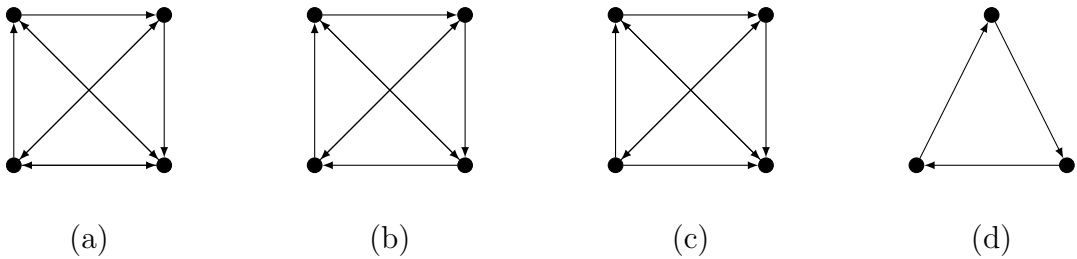
\begin{figure}[htb]
                \center
                \begin{tikzpicture}[>=latex]
                \begin{pgfonlayer}{nodelayer}
                \node [style=blackvertex] (1) at (0,0) {};
                \node [style=blackvertex] (2) at (2,0) {};
                \node [style=blackvertex] (3) at (2,-2) {};
                \node [style=blackvertex] (4) at (0,-2) {};

                \draw[style=arc] (1) to (2);
                \draw[style=arc] (2) to (3);
                \draw[style=arc] (3) to (4);
                \draw[style=arc] (4) to (1);
                \draw[style=arc] (4) to (3);
                \draw[style=arc] (1) to (3);
                \draw[style=arc] (3) to (1);
                \draw[style=arc] (2) to (4);
                \draw[style=arc] (4) to (2);
                \node [style=textbox] at (1, -3) {(a)};

                \node [style=blackvertex] (1) at (4,0) {};
                \node [style=blackvertex] (2) at (6,0) {};
                \node [style=blackvertex] (3) at (6,-2) {};
                \node [style=blackvertex] (4) at (4,-2) {};

                \draw[style=arc] (1) to (2);
                \draw[style=arc] (2) to (3);
                \draw[style=arc] (3) to (4);
                \draw[style=arc] (4) to (1);
                \draw[style=arc] (1) to (3);
                \draw[style=arc] (3) to (1);
                \draw[style=arc] (2) to (4);
                \draw[style=arc] (4) to (2);
                \node [style=textbox] at (5, -3) {(b)};

                \node [style=blackvertex] (1) at (8,0) {};
                \node [style=blackvertex] (2) at (10,0) {};
                \node [style=blackvertex] (3) at (10,-2) {};
                \node [style=blackvertex] (4) at (8,-2) {};

                \draw[style=arc] (1) to (2);
                \draw[style=arc] (2) to (3);
                \draw[style=arc] (4) to (3);
                \draw[style=arc] (4) to (1);
                \draw[style=arc] (1) to (3);
                \draw[style=arc] (3) to (1);
                \draw[style=arc] (2) to (4);
                \draw[style=arc] (4) to (2);
                \node [style=textbox] at (9, -3) {(c)};

                \node [style=blackvertex] (1) at (13,0) {};
                \node [style=blackvertex] (2) at (14,-2) {};
                \node [style=blackvertex] (3) at (12,-2) {};

                \draw[style=arc] (1) to (2);
                \draw[style=arc] (2) to (3);
                \draw[style=arc] (3) to (1);

                \node [style=textbox] at (13, -3) {(d)};

                \end{pgfonlayer}
                \end{tikzpicture}
\caption{\label{Forbidden1}Semicomplete digraphs which are not chordal}
        \end{figure}
\end{center}

A digraph $D$ is called {\em locally semicomplete} if for every vertex $v$, 
$N^-(v)$ and $N^+(v)$ each induces a semicomplete subdigraph in $D$. 
Locally semicomplete digraphs are a popular generalization of semicomplete 
digraphs and have been extensively studied, cf. \cite{lsd,bggv,jbjbook,huang}.
Many properties for semicomplete digraphs hold for locally semicomplete digraphs, 
cf. \cite{lsd}.
However, there are locally semicomplete digraphs which are neither semicomplete 
nor chordal. Any directed cycle consisting of non-symmetric arcs is locally 
semicomplete but not chordal, and is not semicomplete if it has has four or more 
vertices. We will prove that directed cycles with four or more vertices consisting 
of non-symmetric arcs are the only minimal locally semicomplete digraphs which are 
not chordal and which are not semicomplete.

Quasi-transitive digraphs are another well-studied class of digraphs generalizing
semicomplete digraphs, cf. \cite{bh,bh_king,msss,galeana}.
A digraph $D = (V,A)$ is called {\em quasi-transitive} if for any three vertices 
$u, v, w$, $uv \in A$ and $vw \in A$ imply $uw \in A$ or $wu \in A$ (or both), cf.
\cite{bh}. The class of quasi-transitive digraphs contains all {\em transitive 
oriented graphs}. These are the oriented graphs which satisfy the property that
for any three vertices $u, v, w$, $uv \in A$ and $vw \in A$ imply $uw \in A$. 
Equivalently, they are the oriented graphs in which every vertex is a di-simplicial
vertex. Quasi-transitive chordal digraphs are studied recently in \cite{ye}, where  
it is proved that they have the same forbidden subdigraphs as for semicomplete 
chordal digraphs as stated in Theorem \ref{basic}.   

We introduce a new class of digraphs as a common generalization of several 
classes of digraphs including quasi-transitive digraphs and symmetric digraphs.
Since graphs can be viewed as symmetric digraphs, the new class of digraphs 
contains all graphs.

Let $v$ be a vertex and $u, w$ be neighbours of $v$ in a digraph $D$. Then $u, w$ 
are called {\em synchronous} neighbours of $v$ if $u, w$ are both in 
$N^-(v) \setminus N^+(v)$, or in $N^+(v) \setminus N^-(v)$, or in 
$N^-(v) \cap N^+(v)$; otherwise they are called {\em asynchronous} neighbours of 
$v$. We call a digraph $D$ {\em weakly quasi-transitive} if for each vertex $v$ 
of $D$, any two asynchronous neighbours of $v$ are adjacent.    

\begin{center}
        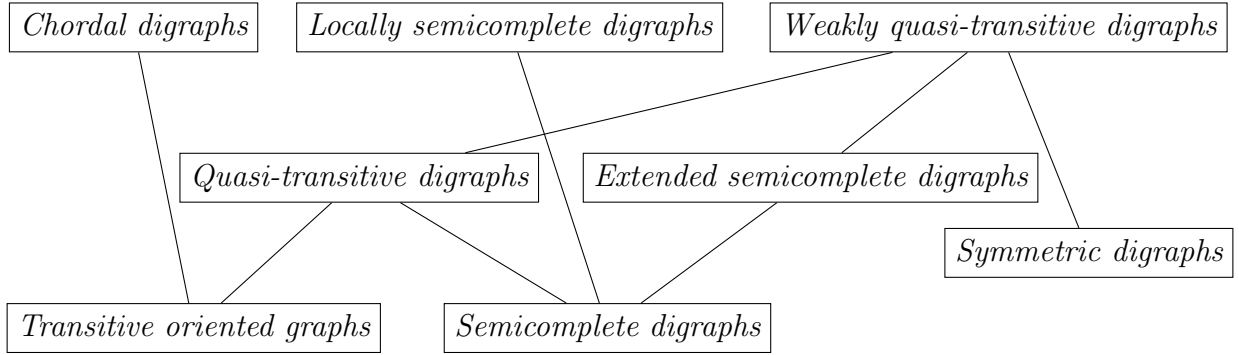
\begin{figure}[h]
                \center
                \begin{tikzpicture}[>=latex]
                \begin{pgfonlayer}{nodelayer}

\node[style=textbox,draw] (2) at (-6,4){{\em Chordal digraphs}};
\node[style=textbox,draw] (3) at (-1,4){{\em Locally semicomplete digraphs}};
\node[style=textbox,draw] (4) at (5.5,4){{\em Weakly quasi-transitive digraphs}};
\node[style=textbox,draw] (5) at (-3,2){{\em Quasi-transitive digraphs}};
\node[style=textbox,draw] (6) at (3,2){{\em Extended semicomplete digraphs}};
\node[style=textbox,draw] (7) at (6.7,1){{\em Symmetric digraphs}};
\node[style=textbox,draw] (8) at (-5.2,0){{\em Transitive oriented graphs}};
\node[style=textbox,draw] (9) at (.3,0){{\em Semicomplete digraphs}};
\draw[style=edge] (2) -- (8) -- (5);
\draw[style=edge] (3) -- (9);
\draw[style=edge] (4) -- (5) -- (9);
\draw[style=edge] (9) -- (6) -- (4);
\draw[style=edge] (7) -- (4);
                \end{pgfonlayer}
                \end{tikzpicture}
\caption{\label{diagram}A containment hierarchy}
        \end{figure}
\end{center}

The class of weakly quasi-transitive digraphs contains all quasi-transitive digraphs
(and hence contains all semicomplete digraphs as well as all transitive oriented 
digraphs). Indeed, suppose $D$ is not weakly quasi-transitive. Then some vertex $v$
has two non-adjacent asynchronous neighbours $u, w$. Since $u, w$ are asynchronous 
neighbours of $v$, one of of $u, w$ is in $N^-(v)$ and the other is in $N^+(v)$. 
Hence $D$ is not a quasi-transitive digraph. Clearly every symmetric digraph is
weakly quasi-transitive. Symmetric digraphs have the property that the neighbours
of each vertex are synchronous and any digraph having this property is weakly
quasi-transitive.  
If a digraph $D$ is weakly quasi-transitive then any digraph obtained from $D$ by 
substituting an independent set for each vertex of $D$ is also weakly 
quasi-transitive. {\em Extended semicomplete digraphs} are the digraphs obtained 
this way from semicomplete digraphs so they are all weakly quasi-transitive. 
Therefore the class of weakly quasi-transitive digraphs simultaneously contains 
quasi-transitive digraphs, symmetric digraphs, and extended semicomplete digraphs. 
Figure~\ref{diagram} depicts a containment hierarchy of the digraph classes 
relevant to this paper.

Let $D$ be a digraph with vertices $v_1, v_2, \dots, v_n$ and let $H_1, H_2, \dots,
H_n$ be vertex-disjoint digraphs. A {\em substitution} of the digraphs $H_i$ for
the vertices $v_i$ in $D$ for each $i$ is a new digraph $D^*$ obtained from 
$H_1, H_2, \dots, H_n$ by adding all possible arcs $xy$ where $x \in V(H_i)$ and 
$y \in V(H_j)$ for each arc $v_iv_j$ in $D$. We use $D[H_1,H_2,\dots,H_n]$ to 
denote the new digraph $D^*$ and also say that it is obtained from $D$ by 
{\em substituting} $H_i$ for $v_i$ for each $i$.
We call $D$ {\em strong} if for any two vertices $x, y$ there is a directed 
$(x,y)$-path and a directed $(y,x)$-path; otherwise we call $D$ {\em non-strong}.

\begin{theorem} \cite{bh} \label{qtd}
Let $D$ be a quasi-transitive digraph. Then the following statements hold:
\begin{enumerate}
\item If $D$ is non-strong, then $D = T[H_1,H_2, \dots, H_n]$ where $T$ is
      a transitive oriented graph and each $H_i$ is a strong quasi-transitive
      digraph.
\item If $D$ is strong, then $D = S[H_1,H_2, \dots, H_n]$ where $S$ is a strong
      semicomplete digraph and each $H_i$ is either a single-vertex digraph
      or a non-strong quasi-transitive digraph.
\qed
\end{enumerate}
\end{theorem}

Thus every quasi-transitive digraph can be obtained from transitive oriented graphs
and semicomplete digraphs recursively by substitutions. Weakly quasi-transitive 
digraphs admit a similar construction. We will show (see Theorem \ref{wqtd}) that 
weakly quasi-transitive digraphs can be constructed recursively
from transitive oriented graphs, symmetric digraphs, and semicomplete digraphs by
substitutions. As a by-product of this recursive construction, we prove 
that the forbidden subdigraphs for weakly quasi-transitive chordal digraphs are
exactly those for semicomplete chordal digraphs. The forbidden subdigraph 
characterization of weakly quasi-transitive chordal digraphs generalizes not only 
the results of \cite{ye} on quasi-transitive chordal digraphs and extended 
semicomplete chordal digraphs but also the classical results on chordal graphs.  


\section{Locally semicomplete chordal digraphs}

Let $D$ be a digraph and $C: v_1v_2 \dots v_kv_1$ be a directed cycle in $D$. 
If there is no arc between $v_i$ and $v_j$ for all $i, j$ with 
$|i-j| \notin \{1,k-1\}$, then the cycle $C$ is called {\em induced} in $D$.
 
\begin{lemma} \label{dicycle} 
If $D$ is a chordal digraph, then $D$ does not contain an induced directed cycle
consisting of non-symmetric arcs and $S(D)$ does not contain an induced directed 
cycle of length $\geq 4$.
\end{lemma}
\pf Suppose that $C$ is either an induced directed cycle in $D$ consisting of 
non-symmetric arcs or an induced directed cycle of length $\geq 4$ in $S(D)$. Then 
the subdigraph of $D$ induced by the vertices of $C$ has no di-simplicial vertex 
and hence is not a chordal digraph. Therefore $D$ is not a chordal digraph.  
\qed

When $S(D)$ contains no induced directed cycle of length $\geq 4$, $S(D)$ is 
a chordal digraph and hence has di-simplicial vertices. The di-simplicial vertices 
of $S(D)$ necessarily contain the di-simplicial vertices of $D$, as observed in
\cite{meister}.  

\begin{lemma}\label{S(D)IsChordal}\cite{meister}
Every di-simplicial vertex of a digraph $D$ is a di-simplicial vertex of $S(D)$.
\qed
\end{lemma}

Suppose that $v$ is a di-simplicial vertex of $S(D)$ but not a di-simplicial vertex
of $D$. Since $v$ is not a di-simplicial vertex of $D$, there exist $u \in N^-(v)$ 
and $w \in N^+(v)$ such that $uw$ is not an arc of $D$ and we shall call such an 
ordered triple $(u,v,w)$ of vertices a {\em violating triple} for $v$. We remark
that a violating triple $(u,v,w)$ exists only if $v$ is a di-simplicial vertex of
$S(D)$ and it certifies that $v$ is not a di-simplicial vertex of $D$. 
We call $v$ {\em type~1} if for every violating triple $(u,v,w)$, both $uv, vw$ 
are non-symmetric and {\em type~2} otherwise.
The following lemma allows us to streamline the selection of violating triples.

\begin{lemma} \label{canonical}
Let $D$ be a locally semicomplete digraph such that $S(D)$ is chordal and $D$ does 
not contain an induced directed cycle consisting of non-symmetric arcs or any 
digraph in Figure \ref{Forbidden1} as an induced subdigraph. Suppose that $(u,v,w)$ 
is a violating triple. Then the following statements hold:
\begin{enumerate}
\item If $uv$ is a non-symmetric arc then there exists a di-simplicial vertex $u'$ 
of $S(D)$ (possibly $u' = u$) such that $(u',v,w)$ is a violating triple and $u'v$ 
is a non-symmetric arc. 
\item If $vw$ is a non-symmetric arc then there exists a di-simplicial vertex $w'$ 
of $S(D)$ (possibly $w' = w$) such that $(u,v,w')$ is a violating triple and $vw'$ 
is a non-symmetric arc. 
\end{enumerate}
\end{lemma}
\pf The two statements can be obtained from each other by reversing the arcs of $D$.
Thus we only prove the first one. Assume $uv$ is non-symmetric. 
Consider $S(D)-(N^-[w] \cap N^+[w])$ where $N^-[w] = N^-(w) \cup \{w\}$ and 
$N^+[w] = N^+(w) \cup \{w\}$. Since $uw$ is not an arc of $D$, 
$u$ is not a vertex in $N^-[w] \cap N^+[w]$ and hence is a vertex of 
$S(D)-(N^-[w] \cap N^+[w])$. Since $S(D)$ is chordal, $S(D)-(N^-[w] \cap N^+[w])$ 
is also chordal. Thus each component of $S(D)-(N^-[w] \cap N^+[w])$ contains
a di-simplicial vertex. Let $u'$ be a di-simplicial vertex of the component of 
$S(D)-(N^-[w] \cap N^+[w])$ that contains $u$. 
Let $u_1u_2 \dots u_k$ where $u_1 = u$ and $u_k = u'$ be a directed path in 
$S(D)-(N^-[w] \cap N^+[w])$. We prove by induction on $k$ that 
$(u_k,v,w)$ is a violating triple and $u_kv$ is a non-symmetric arc of $D$. 

This is true when $k = 1$. So assume $k > 1$, $(u_{k-1},v,w)$ is a violating
triple and $u_{k-1}v$ is a non-symmetric arc of $D$. 
Suppose that $vw$ is a symmetric arc.
Then there is an arc between $u_{k-1}$ and $w$ as they are both in-neighbours of
$v$. Since $u_{k-1}w$ is not an arc of $D$, $wu_{k-1}$ is a non-symmetric arc. 
Thus both $w$ and $u_k$ are in-neighbours of $u_{k-1}$ so there is an arc between 
them. Since $u_k \notin N^-[w] \cap N^+[w]$, $w$ and $u_k$ are joined by 
a non-symmetric arc. Since $u_k$ and $v$ are out-neighbours of $u_{k-1}$, 
they are adjacent. If $u_k$ and $v$ are joined by symmetric arcs, then $u_k$ and 
$w$ are in $N^-(v) \cap N^+(v)$. Since $v$ is a di-simplicial vertex, 
$u_k$ and $w$ are joined by symmetric arcs, which contradicts the fact 
$u_k \notin (N^-[w] \cap N^+[w])$.  
If $vu_k$ or $u_kw$ is an arc of $D$, then the subdigraph of $D$
indcued by $v, w, u_{k-1},u_k$ is Figure \ref{Forbidden1}(a), (b) or (c),
contradicting to our assumption. Hence $vu_k$ is not an arc (i.e., $u_kv$ is
a non-symmetric arc) and $u_kw$ is not an arc of $D$, that is, 
$(u_k,v,w)$ is a violating triple.  
On the other hand, suppose that $vw$ is a non-symmetric arc. Since $u_{k-1}w$ is
not an arc and $D$ does not contain an induced directed cycle consisting of 
non-symmetric arcs, there is
no arc between $u_{k-1}$ and $w$. This implies there is no arc between $u_k$
and $w$ as otherwise $u_{k-1}, w$ are non-adjacent vertices in $N^-(u_k)$ or in
$N^+(u_k)$, which contradicts that $D$ is locally semicomplete. 
The vertices $u_k, v$ are adjacent because they are out-neighbours of $u_{k-1}$.
Since $w$ is an out-neighbour of $v$ and there is no arc between $u_k$ and $w$,
$u_k$ cannot be an out-neighbour of $v$. Therefore $u_kv$ is a non-symmetric arc
of $D$, which means $(u_k,v,w)$ is a violating triple. 
\qed

We call a violating triple $(u,v,w)$ {\em canonical} if $u$ is a di-simplicial 
vertex of $S(D)$ whenever $uv$ is a non-symmetric arc and $w$ is a di-simplicial 
vertex of $S(D)$ whenever $vw$ is a non-symmetric arc.  
Lemma \ref{canonical} ensures that if there is a violating triple for $v$ then
there exists a canonical violating triple for $v$.
In particular, if $v$ is a type~1 vertex then there is a violating triple $(u,v,w)$
for $v$ such that $u, w$ are both di-simplicial vertices of $S(D)$.

\begin{theorem}\label{LSCD}
A locally semicomplete digraph $D$ is chordal if and only if $S(D)$ is chordal and
$D$ does not contain as an induced subdigraph a directed cycle consisting of 
non-symmetric arcs or a digraph in Figure \ref{Forbidden1}.
\end{theorem}
\pf The necessity follows from Theorems \ref{basic} and Lemma \ref{dicycle}. For
the other direction assume that $S(D)$ is chordal and $D$ contains neither a 
directed cycle consisting of non-symmetric arcs nor a digraph in Figure 
\ref{Forbidden1} as an induced subdigraph. To prove $D$ is chordal it suffices to 
show that $D$ has a di-simplicial vertex. Since 
$S(D)$ is chordal, $S(D)$ has di-simplicial vertices. If some di-simplicial vertex 
of $S(D)$ is a di-simplicial vertex of $D$ then we are done. Hence we also assume 
that none of the di-simplcial vertices of $S(D)$ is a di-simplicial vertex of $D$. 

First suppose that $S(D)$ has di-simplicial vertices of type~1. Let $v$ be such a 
vertex. Then there is a violating triple for $v$ and thus by Lemma~\ref{canonical} 
there is a canonical violating triple $(u,v,w)$ for $v$. Note that $u, w$ are both 
di-simplicial vertices of $S(D)$ and $uw$ is not an arc. Since $D$ contains no 
directed cycles consisting of non-symmetric arcs, $wu$ is not an arc and so $u, w$ 
are not adjacent. We claim that the triple $(u,v,w)$ can be chosen so that $u$ is 
type~1. We prove this by contradiction. So assume $u$ is type~2. Then there is 
a canonical violating triple $(u_1,u,w_1)$ for $u$ such that exactly one of 
$u_1u, uw_1$ is a non-symmetric arc.
Suppose first that $u_1u$ is non-symmetric and $uw_1$ is symmetric. 
Since $u_1, w_1$ are both in-neighbours of $u$ and $D$ is locally semicomplete, 
they are adjacent. But $u_1w_1$ is not an arc so $w_1u_1$ is a non-symmetric arc.
There is no arc between $w_1$ and $w$ as otherwise $w,u$ are in-neighbours or 
out-neighbours of $w_1$, which contradicts the fact that they are not adjacent.
Since $w$ is an out-neighbour of $v$ but not adjacent to $w_1$, $w_1$ cannot be 
an out-neighbour of $v$. But $w_1$ and $v$ are adjacent as they are out-neighbours 
of $u$ so $w_1v$ is a non-symmetric arc. Since $u$ is an out-neighbour of $u_1$ but
not adjacent to $w$, $w$ cannot be an out-neighbour of $u_1$. Similarly, $w_1$ is 
an in-neighbour of $u_1$ but not adjacent to $w$, $w$ cannot be an in-neighbour of 
$u_1$. Hence $w$ is not adjacent to $u_1$. There must be an arc between $v$ and 
$u_1$ as they are out-neighbours of $w_1$. But $u_1$ cannot be an out-neighbour of 
$v$ because it is not adjacent to $w$ which is an out-neighbour of $v$. Hence 
$u_1v$ is a non-symmetric arc and $(u_1,v,w)$ is a violating triple. 
Since $(u_1,u,w_1)$ is a canonical violating triple and $u_1u$ is non-symmetric,
$u_1$ is a di-simplicial vertex of $S(D)$ and hence $(u_1,v,w)$ is a canonical
violating triple. A similar proof shows that if $u_1u$ is symmetric and $uw_1$ is 
non-symmetric then $(w_1,v,w)$ is a canonical violating triple. Therefore we have
proved that in the case when $u$ is not type~1 there exists a vertex $x$ (which
is $u_1$ or $w_1$) such that $(x,v,w)$ is a canonical violating triple and the arc 
between $x$ and $u$ is non-symmetric. 
If $x$ is type~1 then it is a desired vertex. Otherwise $x$ is type~2. Repeating
the same argument as above with $x$ replacing $u$ we find the next vertex 
$x'$ which is either a desired vertex or a type~2 vertex such that $(x',v,w)$ is 
a canonical violating triple. Continuing this way in a finite number of steps we 
either find a desired vertex $u$ (i.e., $u$ is type~1 and $(u,v,w)$ is 
a canonical violating triple) or a `circuit' $x_1, x_2, \dots, x_k$, along with
vertices $y_1, y_2, \dots, y_k$, such that for each $i = 1, 2, \dots, k$,

\begin{itemize}
\item $x_i$ is di-simplicial vertex of $S(D)$ of type~2,
\item $(x_i,v,w)$ and $(y_i,v,w)$ are canonical violating triples, and
\item there is a non-symmetric arc between $x_i$ and $x_{i+1}$ and either 
      $(x_{i+1},x_i,y_i)$ or $(y_i,x_i,x_{i+1})$ is a canonical violating triple
      (subscritps are modulo $k$).
\end{itemize} 

Assume the latter occurs and the circuit has the minimum length. Note that 
the vertices $x_1, \dots, x_k, y_1, \dots, y_k$ are in-neighbours of $v$ so they 
are pairwise adjacent. Since $D$ does not contain any digraph in 
Figure \ref{Forbidden1} as an induced subdigraph, the circuit is not a directed 
cycle (consisting of non-symmetric arcs). Hence we may assume without loss of 
generality that $x_1x_2, x_1x_k$ are non-symmetric arcs. 
Then $(y_1,x_1,x_2)$ and $(x_1,x_k,y_k)$ are canonical violating
triples. If $x_2y_k$ is a non-symmetric arc then $x_1,x_2,y_k$ induce 
Figure~\ref{Forbidden1}(d), a contradiction to assumption. If $x_2y_k$ is symmetric
then $x_1,x_2,y_1,y_k$ induce Figure~\ref{Forbidden1}(a), (b) or (c), also a 
contradiction. So $y_kx_2$ is a non-symmetric arc. Since $x_k$ is a di-simplicial 
vertex of $S(D)$ and $y_k$ is a adjacent to $x_k$ but not to $x_2$ in $S(D)$, 
$x_2$ is not adjacent to $x_k$ in $S(D)$, that is, the arc between $x_2$ and $x_k$ 
is non-symmetric. 
If $x_2x_k$ is an arc then $x_2, \dots, x_k$ would be a shorter circuit, 
a contradiction to our choice of circuit. So $x_kx_2$ is a non-symmetric arc. 
There is an arc between $y_1$ and $x_k$ as they are out-neighbours of $x_1$.
If $y_1x_k$ is non-symmetric, then $y_1,x_k,x_2$ induce Figure~\ref{Forbidden1}(d)
and if $x_ky_1$ is non-symmetric, then $x_1,x_k,y_1,y_k$ induce 
Figure~\ref{Forbidden1}(a), (b) or (c), a contradiction to assumption.
Hence $y_1x_k$ is a symmetric arc and $x_2, \dots, x_k$ is a shorter circuit,
which is also a contradiction. Therefore for every type~1 vertex $v$ there exists
a canonical violating triple $(u,v,w)$ such that $u$ is a type~1 vertex.
This implies that there exists a directed cycle on type~1 vertices consisting of 
non-symmetric arcs. Assume that $v_1, v_2, \dots, v_t$ is the shortest such cycle.
Since $D$ does not contain an induced directed cycle consisting of non-symmetric
arcs, $t > 3$ and there is a symmetric arc joining a pair of non-consecutive 
vertices of the cycle. 
Without loss of generality assume $v_1v_s$ is a symmetric 
arc of the shortest distance along the cycle, that is, $v_i, v_j$ are not adjacent 
for all $1 \leq i < j-1 \leq s-1$ except $i=1$ and $j = s$. Since $v_2$ and $v_s$ 
are out-neighbour of $v_1$, they are adjacent. This implies that $s = 3$ and so
$v_1v_3$ is a symmetric arc. Hence $(v_2,v_3,v_1)$ is a violating triple in which 
$v_3v_1$ is a symmetric arc, which contradicts the assumption that $v_3$ is 
a type~1 vertex. Therefore $S(D)$ has no type~1 di-simplicial vertex, that is,
every di-simplicial vertex of $S(D)$ is type~2.

Let $v$ be a di-simplicial vertex of $S(D)$. Since $v$ is type~2, there is a
canonical violating triple $(u,v,w)$ such that exactly one of $uv, vw$ is a 
non-symmetric arc. If $uv$ is non-symmetric then $u$ is a di-simplicial vertex
of $S(D)$. If $vw$ is non-symmetric then $w$ is a di-simplicial. This implies
that for each di-simplicial vertex of $S(D)$ there is a di-simplicial vertex 
$z$ of $S(D)$ such that $z, v$ are part of a canonical violating triple for $v$ 
and the arc between $v$ and $z$ is non-symmetric. It follows that there exists
a `circuit' $z_1, z_2, \dots, z_r$, along with vertices $w_1, w_2, \dots, w_r$, 
such that for each $i = 1, 2, \dots, r$,

\begin{itemize}
\item $z_i$ is a di-simplicial vertex of $S(D)$ of type~2,
\item either $(z_{i+1},z_i,w_i)$ or $(w_i,z_i,z_{i+1})$ is a canonical violating
      triple,
\item the arc between $z_i$ and $z_{i+1}$ is non-symmetric and the arcs between
      $w_i$ and $z_i$ are symmetric (subscripts are modulo $r$).
\end{itemize} 

We again assume that the circuit is chosen to have the minimum length. 
Suppose $r = 2$. If the non-symmetric arc between $z_1$ and $z_2$ is $z_1z_2$,
then $(w_1,z_1,z_2)$ and $(z_1,z_2,w_2)$ are the canonical violating triples where
$w_1z_1$ and $z_2w_2$ are symmetric arcs. Neither $w_1z_2$ nor $z_1w_2$ is an
arc. Since $w_1$ and $z_2$ are out-neighbours of $z_1$, they are adjacent so
$z_2w_1$ is a non-symmetric arc. Similarly, $w_2z_1$ is a non-symmetric arc. 
There is an arc between $w_1$ and $w_2$ as they are in-neighbours of $z_1$.
Depending the arcs between $w_1$ and $w_2$, the subdigraph induced by 
$z_1, z_2, w_1, w_2$ is Figure~\ref{Forbidden1}(a), (b) or (c), which contradicts
the assumption. The same conclusion holds if the non-symmetric arc between 
$z_1$ and $z_2$ is $z_2z_1$. So $r \geq 3$. 

Suppose that $z_1z_2 \dots z_rz_1$ is a directed cycle. Since $D$ does not contain 
an induced directed cycle consisting of non-symmetric arcs, $r > 3$ and there is 
a symmetric arc between a pair of non-consecutive vertices of the cycle. 
Without loss of generality assume $z_1z_s$ is a symmetric arc of the shortest 
distance along the cycle, that is, $z_i, z_j$ are not adjacent for all 
$1 \leq i < j-1 \leq s-1$ except $i=1$ and $j = s$. 
Since $z_2$ and $z_s$ are out-neighbour of $z_1$, they are adjacent. 
This implies that $s = 3$ and so $z_1z_3$ is a symmetric arc. Since $z_3$ and $z_r$
are in-neighbours of $z_1$, they are adjacent. Since $z_3$ and 
$w_r$ are out-neighbours of $z_1$, they are adjacent. The arcs between $z_3$ and
$w_r$ cannot be symmetric as otherwise $z_1$ and $w_r$ are both neighbours of 
$z_3$ in $S(D)$ but $z_1w_r$ is a non-symmetric arc, which contradicts the fact
that $z_3$ is a di-simplicial vertex of $S(D)$. So $z_3$ and $w_r$ are joined by
a non-symmetric arc. If $w_rz_3$ is a non-symmetric arc, then the subdigraph 
induced by $z_1, z_3, z_r, w_r$ is Figure~\ref{Forbidden1}(a), (b) or (c), 
a contradiction. Hence $z_3w_r$ is a non-symmetric arc. The arc between $z_3$ and 
$z_r$ must be non-symmetric as otherwise $z_3$ and $w_r$ are non-adjacent 
neighbours of $z_r$ in $S(D)$, which contradicts the fact that $z_r$ is 
a di-simplicial vertex of $S(D)$. 
If $z_3z_r$ is a non-symmetric arc then $z_1, z_3, z_r, w_r$ induce 
Figure~\ref{Forbidden1}(c), a contradiction. On the other hand, if $z_rz_3$ is a 
non-symmetric arc, then $z_3, \dots, z_r$ would be a directed cycle of 
length shorter than $r$ consisting of non-symmetric arcs, which contradicts 
the choice of circuit. Therefore $z_1z_2 \dots z_rz_1$ is not a directed cycle.
Hence we may assume without loss of generality that $z_1z_2$ and $z_1z_r$ are 
non-symmetric arcs.   

Since $z_1z_2$ and $z_1z_r$ are non-symmetric arcs, $(w_1,z_1,z_2)$ and 
$(z_1,z_r,w_r)$ are canonical violating triples. Since $z_2$ and $z_r$ are 
out-neighbours of $z_1$, they are adjacent. So $z_2$ is an in-neighbour or 
an out-neighbour of $z_r$. Combining this with the fact that $w_r$ is both
an in-neighbour and an out-neighbour of $z_r$ we see that $z_2$ and $w_r$ are
adjacent. If $z_2$ and $w_r$ are joined by symmetric arcs then $z_1,z_2,w_1,w_r$
induced Figure~\ref{Forbidden1}(a), (b) or (c), a contradiction. So $z_2$ and $w_r$ are joined by a non-symmetric arc. If $z_2w_r$ is a non-symmetric arc, then
$z_1,z_2,w_r$ induce Figure~\ref{Forbidden1}(d), a contradiction. Hence 
$w_rz_2$ is a non-symmetric arc. This means that $w_r$ is not adjacent to $z_2$
in $S(D)$. However, $w_r$ is adjacent to $z_r$ in $S(D)$ and $z_r$ is a 
di-simplicial vertex of $S(D)$. It follows that $z_2$ and $z_r$ are joined by
a non-symmetric arc. If $z_2z_r$ is a non-symmetric arc, then $z_2, \dots, z_r$
would be a shorter circuit, a contradiction to our choice. So $z_rz_2$ is 
a non-symmetric arc. Since $w_1$ and $z_r$ are out-neighbours of $z_1$, they
are adjacent. If $w_1$ and $z_r$ are joined by symmetric arcs, then 
again $z_2, \dots, z_r$ would be a shorter circuit, a contradiction. So
$w_1$ and $z_r$ are joined by a non-symmetric arc. It cannot be $w_1z_r$ as
otherwise $w_1,z_r,z_2$ induce Figure~\ref{Forbidden1}(d), a contradiction. 
Hence $z_rw_1$ is a non-symmetric arc. The subdigraph induced by $z_1,z_r,w_1,w_r$
is Figure~\ref{Forbidden1}(a), (b) or (c), a contradiction.    
Therefore, $D$ has a di-simplicial vertex. This completes the proof.	
\qed

\section{Weakly quasi-transitive digraphs}

According to Theorem \ref{qtd}, transitive oriented graphs and semicomplete 
digraphs are basic building blocks for quasi-transitive digraphs. Using these 
blocks one can form a class $\cal Q$ of digraphs as follows:  

\begin{enumerate}
\item Each transitive oriented graph is in $\cal Q$.
\item Each semicomplete digraph is in $\cal Q$.
\item If $D, H_1, H_2, \dots, H_n \in \cal Q$, then $D[H_1,H_2,\dots,H_n] \in 
      \cal Q$, provided that $H_i$ is a single-vertex digraph when the vertex $v_i$ 
      for which $H_i$ is substituted is incident with a symmetric arc for each $i$.
\end{enumerate}

Transitive oriented graphs and semicomplete digraphs are quasi-transitive. 
Moreover, the substitution operation for defining $\cal Q$ maintain the property 
of being quasi-transitive. Hence the digraphs in $\cal Q$ are all quasi-transitive.
Theorem \ref{qtd} ensures that every quasi-transitive digraph can be obtained from 
transitive oriented graphs and semicomplete digraphs by substitutions. Therefore
we have the following:    

\begin{corollary} \label{Q}
The class $\cal Q$ consists of quasi-transitive digraphs.
\qed
\end{corollary}

Interestingly, weakly quasi-transitive digraphs can also be constructed in a similar
way from transitive oriented graphs, semicomplete digraphs and symmetric digraphs. 

Let $\cal W$ be the class of digraphs defined as follows:

\begin{enumerate}
\item each transitive oriented graph is in $\cal W$;
\item each semicomplete digraph is in $\cal W$;
\item each symmetric digraph is in $\cal W$;
\item if $D$ is in $\cal W$ then any digraph obtained from $D$ by substituting
      digraphs of $\cal W$ for the vertices of $D$ is in $\cal W$.
\end{enumerate}

A {\em module} in a digraph $D$ is an induced subgraph $H$ of $D$ such that for 
any vertex $x$ not in $H$, either $x$ is adjacent to no vertex in $H$ or 
the vertices in $H$ are synchronous neighbours of $x$. A module is called 
{\em trivial} if it has only one vertex or is the entire digraph $D$ and 
{\em non-trivial} otherwise. An {\em oriented path} in $D$ is a sequence of 
vertices $v_1, v_2, \dots, v_k$ such that $v_i$ and $v_{i+1}$ are joined by 
a non-symmetric arc for each $i = 1, 2, \dots, k-1$. 

\begin{theorem} \label{wqtd}
The class $\cal W$ consists of weakly quasi-transitive digraphs.
\end{theorem}
\pf Transitive oriented graphs and semicomplete digraphs are quasi-transitive, so
they are weakly quasi-transitive. Symmetric digraphs are also weakly 
quasi-transitive because any vertex in a symmetric digraph has only synchronous 
neighbours. To prove the rest of digraphs in $\cal W$ are all weakly 
quasi-transitive, let $D^* = D[H_1,H_2,\dots,H_n]$ where $D,H_1,H_2, \dots,H_n$ are
weakly quasi-transitive. Consider three vertices $u,v,w$ where $u,w$ are 
asynchronous neighbours of $v$. Assume $u \in V(H_i)$, $v \in V(H_j)$ and 
$w \in V(H_k)$. 
If $i=j=k$ then $u,w$ are adjacent as $H_i$ is weakly quasi-transitive.  
Suppose $i= j \neq k$. Since $v$ and $w$ are adjacent, each vertex of $H_i$ is
adjacent to all vertices of $H_k$ and in particular, $u$ is adjacent to $w$. 
Similarly, if $i \neq j = k$, then $u$ and $w$ are adjacent. Suppose that 
$i \neq j \neq k$. Then $i \neq k$ because $u$ and $w$ are asynchronous neighbours
of $v$. Since $D$ is weakly quasi-transitive, the two vertices of $D$ corresponding
to $H_i$ and $H_k$ are adjacent so $u$ and $w$ are adjacent. Hence all digraphs
in $\cal W$ are weakly quasi-transitive.

We prove by induction on number of vertices that every weakly quasi-transitive 
digraph is in $\cal W$. Let $D$ be a weakly quasi-transitive with $n$
vertices. Assume that every weakly quasi-transitive digraph with fewer than $n$
vertices is in $\cal W$. If $D$ is quasi-transitive or symmetric then it is in 
$\cal W$. So assume that $D$ is neither quasi-transitive nor symmetric. Since $D$ 
is not quasi-transitive, there exist vertices $u, v, w$ with $u \in N^-(v)$ and 
$w \in N^+(v)$ such that $u$ and $w$ are not adjacent in $D$. Thus $u$ and $w$ are
non-adjacent neighbours of $v$. Since $D$ is weakly quasi-transitive, any 
two asynchronous neighbours of $v$ are adjacent. Hence $u$ and $w$ are synchronous
neighbours of $v$, which implies $u$ and $w$ are both in $N^+(v) \cap N^-(v)$. 

Suppose $H$ is a non-trivial module in $D$. Let $D'$ be the digraph obtained from 
$D$ by deleting all vertices of $H$ except one. Then $D = D'[H_1,H_2,\dots,H_k]$ 
where $H_1 = H$ and each $H_i$ with $i \geq 2$ is a single-vertex digraph. 
The digraphs $D', H_1, \dots, H_k$ each has fewer than $n$ vertices and is weakly 
quasi-transitive and hence they are in $\cal W$. This means that $D$ is obtained 
from digraphs in $\cal W$ by substitution and by definition $D$ is in $\cal W$. 
Thus, it suffices to show that there is a non-trivial module in $D$.

Let $R$ be the subdigraph of $D$ induced by $N^+(v) \cap N^-(v)$. Then $u$ and $w$ 
are a pair of non-adjacent vertices in $R$.  Let $M_1$ be the subdigraph of $R$ 
induced by the vertices which are connected to $u$ by paths in $\overline{U(R)}$.
Clearly, $M_1$ contains $u$ and $w$ but not $v$. Suppose $x$ is a vertex in 
$N^+[v] \cup N^-[v]$ but not in $M_1$. We claim that $x$ is completely adjacent 
to $M_1$. Indeed, if $x \in N^+[v] \cap N^-[v]$, then the definition of $M_1$ 
implies that $x$ is completely adjacent to $M_1$. On the other hand, if 
$x \in N^+(v) \oplus N^-(v)$, then $x$ and any vertex of $M_1$ are asynchronous 
neighbours of $v$ so $x$ is also completely adjacent to $M_1$.  By definition any 
two vertices of $M_1$ are connected by a path in $\overline{U(M_1)}$. In such 
a path any two consecutive vertices are not adjacent in $D$ and hence are 
synchronous neighbours of $x$. It follows that the vertice of $M_1$ are 
synchronous neighbours of $x$. Suppose $x \notin N^+[v] \cup N^-[v]$. If $x$ is
adjacent to some vertex $y$ in $M_1$, then $x$ and $v$ are non-adjacent neighbours 
of $y$ and hence they must be synchronous neighbours of $y$. The fact that 
$v$ is joined to $y$ by symmetric arcs implies $x$ is joined to $y$ by symmetric
arcs. Thus if $x$ is completely adjacent to $M_1$ then the vertices of $M_1$ are 
synchronous neighbours of $x$. It follows that $M_1$ is a module if for each 
$x \notin N^+[v] \cup N^-[v]$, either $x$ is adjacent to no vertex in $M_1$ or 
completely adjacent to $M_1$. We may assume $M_1$ is not a module as otherwise we 
are done. This means that there exist vertices $x, y, y'$ with 
$x \notin N^+[v] \cup N^-[v]$ and $y, y' \in M_1$ such that $x$ is adjacent to 
$y$ but not to $y'$. These three vertices $x, y, y'$ along with $M_1$ will be 
refered to in the rest of proof.

Suppose $N^+(v) \oplus N^-(v) \neq \emptyset$. Any vertex in $N^+(v) \oplus N^-(v)$
is a neighbour of $v$ asynchronous to those of $v$ in $N^+(v) \cap N^-(v)$. Hence
every vertex in $N^+(v) \oplus N^-(v)$ is completely adjacent to 
$N^+(v) \cap N^-(v)$ and in particular to $M_1$. 
Suppose that the arcs between $N^+(v) \oplus N^-(v)$ and $M_1$ are all symmetric. 
Let $M_2$ be the subdigraph of $D$ induced by vertices 
which are connected to $v$ by oriented paths. Clearly, $M_2$ contains $v$ and 
all vertices in $N^+(v) \oplus N^-(v)$. We show that $x$ is not a vertex in $M_2$.
Suppose not; there is an oriented path connecting $x$ and $v$. Then there must
exists an oriented path connecting $x$ and a vertex in $N^+(v) \oplus N^-(v)$. 
Let $a_1\sim a_2 \sim \dots \sim a_s$ be such a path where $a_1 = x$ and
$a_s \in N^+(v) \oplus N^-(v)$. Note that $a_s$ is joined to each vertex of 
$M_1$ by symmetric arcs and $a_1$ ($=x$) is not adjacent to $y'$ (in $M_1$). 
Let $j$ be the largest subscript such that $a_j$ is not adjacent to some vertex 
$y''$ of $M_1$. Then $j < k$ and $a_j \notin N^+[v] \cup N^-[v]$. Since $a_j$ and 
$a_{j+1}$ are joined by a non-symmetric arc, $a_{j+1} \notin N^+[v] \cap N^-[v]$.
Either $a_{j+1} \in N^+(v) \oplus N^-(v)$ or $a_{j+1} \notin N^+[v] \cup N^-[v]$.
In either case $a_{j+1}$ is joined to each vertex of $M_1$ by symmetric arcs. 
Thus $a_j$ and $y''$ are non-adjacent asynchronous neighbours of $a_{j+1}$,
contradicting the assumption that $D$ is weakly quasi-transitive. So $x$ is not 
a vertex of $M_2$. We show that $M_2$ is a module. Let $z$ be a vertex not in $M_2$.
By definition $z$ cannot be joined to any vertex of $M_2$ by a non-symmetric arc. 
Suppose $z$ is joined to some vertex $h$ of $M_2$ by symmetric arcs. Since $h$ can 
reach every other vertex of $M_2$ by an oriented path, following such a path we see
that $z$ is joined to every vertex in the path by symmetric arcs. 
Hence the vertices of $M_2$ are synchronous neighbours of $z$. Therefore $M_2$ is 
a non-trivial module in $D$.   

Suppose now that the arcs between $N^+(v) \oplus N^-(v)$ and $M_1$ are not all 
symmetric. Let $M_3$ be a subdigraph of $D$ induced by the vertices defined
recursively as follows:   
\begin{itemize}
\item $u$ is a vertex in $M_3$;
\item if $h$ is a vertex in $N^+(v) \cap N^-(v)$ that is not adjacent to a vertex
      in $M_3$ then $h$ is a vertex in $M_3$;
\item if $h$ is not in $N^+(v) \cap N^-(v)$ that is joined to a vertex in $M_3$ 
      by symmetric arcs then $h$ is a vertex in $M_3$.
\end{itemize}
\noindent It is easy to see that $M_3$ contains $u,v,w,x$ and all vertices of $M_1$.
Let $b$ be a vertex in $N^+(v) \oplus N^-(v)$ which is joined to a vertex in $M_1$
by a non-symmetric arc. Assume that $b \in N^-(v) \setminus N^+(v)$. From the above
we know that the vertices of $M_1$ are synchronous neighbours of $b$. In particular,
$y, y'$ are synchronous neighbours of $b$. The vertex $y$ is joined to $b$ by
a non-symmetric arc and joined to $x$ by symmetric arcs. Thus $b$ and $x$ are
asynchronous neighbours of $y$ and hence they must be adjacent. 
So $x$ and $v$ are neighbours of $b$. Since $x$ and $v$ are
not adjacent, they are synchronous neighbours of $b$. 
Since $b \in N^-(v) \setminus N^+(v)$, $bv$ is a non-symmetric arc, so $bx$ is 
also a non-symmetric arc. Since $bx$ is a non-symmetric arc and $x, y'$ are 
non-adjacent neighbours of $b$, $by'$ is also a non-symmetric arc. 
The fact that the vertices of $M_1$ are synchronous neighbours of $b$ so
there is a non-symmetric arc from $b$ to every vertex in $M_1$. 
Similarly, if $b \in N^+(v) \setminus N^-(v)$ is joined to a vertex in $M_1$ by 
a non-symmetric arc then $xb$ is a non-symmetric arc and there is a non-symmetric 
arc from every vertex of $M_1$ to $b$.  
 
We claim that $b$ is not a vertex in $M_3$. Suppose not; $b$ is in $M_3$. 
By the definition of $M_3$ there exists a sequence of vertices 
$h_0, h_1, \dots, h_t$ where $h_0 = y$ and $h_t = b$ such that for each
$i > 0$, $h_i \in N^+(v) \cap N^-(v)$ implies that $h_i$ is not adjacent to 
$h_{i-1}$, and $h_i \notin N^+(v) \cap N^-(v)$ implies $h_i$ is joined to $h_{i-1}$
by symmetric arcs. We choose such a vertex $b$ so that the sequence is as short
as possible. Assume $b \in N^-(v) \setminus N^+(v)$. Then $b$ ($=h_t$) is joined to 
$h_{t-1}$ by symmetric arcs. We claim $h_{t-1} \in N^+(v) \cap N^-(v)$.       
Indeed, since $b$ is joined to $h_{t-1}$ by symmetric arcs, 
$h_{t-1} \in N^+(v) \cup N^-(v)$. Suppose $h_{t-1} \in N^-(v) \setminus N^+(v)$.
The choice of $b$ implies that there can only be symmetric arcs between $h_{t-1}$
and $M_1$. Since $h_{t-1}$ and $x$ are asynchronous neighbours of $b$, they are 
adjacent. In particular, $h_{t-1}x$ is a non-symmetric arc. Thus
$x, y'$ are non-adjacent asynchronous neighbours of $h_{t-1}$, a contradiction.
So $h_{t-1} \notin N^-(v) \setminus N^+(v)$. A similar proof shows 
$h_{t-1} \notin N^+(v) \setminus N^-(v)$. So $h_{t-1} \in N^+(v) \cap N^-(v)$. 
Since $b$ is joined to $h_{t-1}$ by symmetric arcs and joined to each vertex of 
$M_1$ by a non-symmetic arc, $h_{t-1} \notin M_1$ and thus $t > 2$. Hence
$h_{t-1}$ is not adjacent to $h_{t-2}$ and is completely adjacent to $M_1$. 
If $h_{t-2} \in N^+[v] \cup N^-[v]$, then $h_{t-2}$ must be in $N^+(v) \cup N^-(v)$
and hence adjacent to $b$. Thus $h_{t-1}, h_{t-2}$ are neighbours of $b$.
Since $h_{t-1}, h_{t-2}$ are not adjacent, they are synchronous neighbours of $b$,
which implies $b$ is joined to $h_{t-2}$ by symmetric arcs. This contradicts 
the choice of the sequence as $h_0, h_1, \dots, h_{t-2}, b$ is a shorter 
sequence. So $h_{t-2} \notin N^+[v] \cup N^-[v]$. Let $\ell$ be the largest
integer such that $h_{t-2}, \dots, h_{t-\ell}$ are not in $ N^+[v] \cup N^-[v]$.  
Then $h_{t-i}$ is joined to $h_{t-i-1}$ by symmetric arcs for each 
$i=2, \dots, \ell$. We must have $h_{t-\ell-1} \in N^+(v) \cap N^-(v)$.
The vertex $b$ is not adjacent to $h_{t-2}$ as otherwise $h_{t-1}, h_{t-2}$ are
non-adjacent asynchronous neighbours of $b$, a contradiction. 
For the same reason, we see that $b$ is not adjacent to $h_{t-i}$ for each
$i = 2, \dots, \ell$. Since $b, h_{t-\ell-1}$ are asynchronous neighbours of $v$,
they are adjacent. They must be joined by symmetric arcs, as otherwise 
$b, h_{t-\ell}$ are non-adjacent asynchronous neighbours of $h_{t-\ell-1}$,
a contradiction. But this contradicts the choice of the sequence because 
$h_0, h_1, \dots, h_{t-\ell-1}, b$ is a shorter sequence. Therefore $b$ is not
a vertex in $M_3$. So if $b \in N^-(v) \setminus N^+(v)$ is joined to a vertex in
$M_1$ with a non-symmetric arc then $b \notin M_3$ and there is a non-symmetric arc
from $b$ to every vertex in $M_1$. A similar proof shows that if 
$b \in N^+(v) \setminus N^-(v)$ is joined to a vertex in $M_1$ with a non-symmetric
arc then $b \notin M_3$ and there is a non-symmetric arc from each vertex of $M_1$ 
to $b$.

We show that $M_3$ is a module. Let $z$ be a vertex that is not in $M_3$.
For each vertex $h \in M_3$, there is a sequence of vertices $h_0,h_1, \dots, h_t$ 
where $h_0 = y$ and $h_t = h$ such that for each $i > 0$, 
$h_i \in N^+(v) \cap N^-(v)$ implies that $h_i$ is not adjacent to 
$h_{i-1}$, and $h_i \notin N^+(v) \cap N^-(v)$ implies $h_i$ is joined to 
$h_{i-1}$ by symmetric arcs. 
Suppose first that $z \in N^-(v) \setminus N^+(v)$. We know from 
the above that $zx$ is a non-symmetric arc and $zh$ is a non-symmetric arc for all 
$h \in M_1$. In particular, $zy$ ($=zh_0$) is a non-symmetric arc.  
Suppose $k > 0$ and $zh_{k-1}$ is a non-symmetric arc. 
If $h_k \in N^+(v) \cap N^-(v)$, then $h_{k-1}, h_k$ are non-adjacent neighbours of
$z$ so $zh_k$ is a non-symmetric arc.
If $h_k \notin N^+(v) \cap N^-(v)$, then $z, h_k$ are asynchronous neighbours of 
$h_{k-1}$ so they are adjacent. There are two cases. Either 
$h_k \in N^+(v) \oplus N^-(v)$ or $h_k \notin N^+(v) \cup N^-(v)$. 
If $h_k \notin N^+(v) \cup N^-(v)$, then clearly $zh_k$ is a non-symmetric arc.  
Assume $h_k \in N^+(v) \oplus N^-(v)$. Since $h_k \notin M_3$, $h_k$ is joined to
each vertex in $M_1$ by symmetric arcs. In particular, $h_k$ is joined to $y'$ by
symmetric arcs. Since $y'$ is not adjacent to $x$, $h_k$ and $x$ cannot be adjacent
as otherwise $y'$ and $x$ are non-adjacent asynchronous neighbours of $h_k$,
a contradiction. Hence $h_k$ and $x$ are synchronous neighbours of $z$. 
Since $zx$ is a non-symmetric arc, $zh_k$ is a non-symmetric arc. 
Therefore $zh$ is a non-symmetric arc for all $h \in M_3$. A similar proof shows
that if $z \in N^+(v) \setminus N^-(v)$ then $hz$ is a non-symmetric arc for all
$h \in M_3$.
Suppose next that $z \in N^+(v) \cap N^-(v)$. Since $z$ is not in $M_3$, $z$ is 
adjacent to
every vertex in $M_3$. In particular, $z$ is adjacent to $x$. Note that $z$ and
$x$ are joined by symmetric arcs. Since $x$ and $y'$ are not adjacent, $z$ is
adjacent to $y'$ by symmetric arcs. This implies $z$ is also joined to $y$ by
symmetric arcs. Suppose that $k > 0$ and $z$ is joined to $h_{k-1}$ by symmetric
arcs. If $h_k \notin N^+(v) \cup N^-(v)$ then clearly $z$ is joined to $h_k$ by
symmetric arcs. 
If $h_k \in N^+(v) \cap N^-(v)$, then $h_k$ is not adjacent to $h_{k-1}$ 
and thus $h_k, h_{k-1}$ are non-adjacent neighbours of $z$. Since $z$ is joined 
to $h_{k-1}$ by symmetric arcs, $z$ is joined to $h_k$ by symmetric arcs. 
If $h_k \in N^+(v) \oplus N^-(v)$, then $h_k$ is joined to $y'$ by symmetric 
arcs. Since $y'$ and $x$ are not adjacent, $h_k$ and $x$ are not adjacent. 
Thus $h_k$ and $x$ are non-adjacent neighbours of $z$, which implies $z$ is 
joined to $h_k$ by symmetric arcs.  
Suppose now that $z \notin N^+[v] \cup N^-[v]$. Since $z$ is not in $M_3$, it is 
not adjacent to any vertex in $M_1$. In particular, $z$ is not adjacent to $y$. 
Suppose that $k > 0$ and $z$ is not adjacent to $h_{k-1}$.
If $h_k \in N^+(v) \cap N^-(v)$, then $z$ is not adjacent to $h_k$ as otherwise
$z$ is joined to $h_k$ by symmetric arcs, which implies $z \in M_3$, a 
contradiction to assumption. If $h_k \notin N^+(v) \cap N^-(v)$, then $h_k$
is joined to $h_{k-1}$ by symmetric arcs. Since $z$ is not adjacent to $h_{k-1}$,
$z$ cannot be joined to $h_k$ by a non-symmetric arc. Since $z \notin M_3$ and 
$h_k \in M_3$, $z$ cannot be joined to $h_k$ by symmetric arcs. Hence $z$ is not 
adjacent to $h_k$.

The only case remaining is that $N^+(v) \oplus N^-(v) = \emptyset$. Since $D$ is 
not a symmetric digraph, it has a non-symmetric arc. Suppose $fg$ is a 
non-symmetric arc in $D$. Let $M_4$ be the subdigraph induced by the vertices
which are connected to $f$ by oriented paths. Then any two vertices in $M_4$ are 
connected by an oriented path. Since $N^+(v) \oplus N^-(v) = \emptyset$, 
there is no oriented path connecting $f$ and $v$. So $v$ is not a vertex in $M_4$. 
Suppose $z$ is not in $M_4$ but is adjacent to a vertex $h$ in $M_4$. Then $z$ is 
joined to $h$ by symmetric arcs. Each vertex of $M_4$ is connected to $h$ by 
an oriented path. Following these oriented paths we see that $z$ is joined to each 
vertex of $M_4$ by symmetric arcs and hence the vertices of $M_4$ are synchronous 
neighbours of $z$. Therefore $M_4$ is a non-trivial module. 
\qed

The class of weakly quasi-transitive digraphs strictly contains quasi-transitive
digraphs and extended semicomplete digraphs, which in turn as classes strictly 
contain all semicomplete digraphs. Surprisingly, these four classes of digraphs 
share the same forbidden subdigraphs for being chordal. 

\begin{theorem} \label{ng}
A weakly quasi-transitive digraph $D$ is chordal if and only if $S(D)$ is chordal
and $D$ does not contain any digraph in Figure \ref{Forbidden1} as an induced 
subdigraph.
\end{theorem}
\pf If $D$ is chordal then it does not contain any digraph in Figure 
\ref{Forbidden1} as an induced subdigraph. Suppose $D$ does not contain any
digraph in Figure \ref{Forbidden1} as an induced subdigraph. We prove by induction
on the number of vertices that $D$ is chordal. It suffices to show that $D$ has 
a di-simplicial vertex. This is true if $D$ is a transitive oriented graph, 
a semicomplete digraph, or a symmetric digraph. Assume $D$ is a weakly
quasi-transitive digraph but not a transitive oriented graph, 
a semicomplete digraph, or a symmetric digraph. For the inductive hypothesis,
assume that any induced subdigraph of $D$ with fewer vertices than $D$ has 
a di-simplicial vertex. By Theorem \ref{wqtd}, $D = D'[H_1,H_2,\dots,H_n]$
where $D'$ and one of $H_i$'s have at least two vertices. Then $D'$ and each $H_i$
is an induced subdigraph of $D$ with fewer vertices than $D$ and by the inductive
hypothesis each of them has a di-simplicial vertex. Suppose that $v$ is 
a di-simplicial vertex of $D'$ and $H_j$ is substituted for $v$. Then it is 
easy to verify that a di-simplicial vertex of $H_j$ is a di-simplicial vertex
of $D$. 
\qed

\begin{corollary} \cite{ye}
Let $D$ be a quasi-transitive digraph or an extended semicomplete digraph.
Then $D$ is chordal if and only if $S(D)$ is chordal and $D$ does not contain 
any digraph in Figure \ref{Forbidden1} as an induced subdigraph.
\qed
\end{corollary}

Since graphs can be viewed as symmetric digraphs which are a subclass of the class 
of weakly quasi-transitive digraphs and none of the digraphs in Figure 
\ref{Forbidden1} is symmetric, Theorem \ref{ng} implies that the cycles of
length $\geq 4$ are precisely the forbidden induced subgraphs of chordal graphs.


\begin{thebibliography}{100}

\bibitem{lsd} J. Bang-Jensen, Locally semicomplete digraphs: A generalization
of tournaments, J. Graph Theory 14 (1990) 371 - 390. 

\bibitem{bggv} J. Bang-Jensen, Y. Guo, G. Gutin, and L. Volkmann, 
A classification of locally semicomplete digraphs,
Discrete Math. 167-168 (1997) 101 - 114.

\bibitem{jbjbook} J. Bang-Jensen and G. Gutin, {\em Classes of Directed Graphs},
Springer Monographs in Mathematics (2018).

\bibitem{bh} J. Bang-Jensen and J. Huang, Quasi-transitive digraphs,
J. Graph Theory 20 (1995) 141 - 161. 

\bibitem{bh_king} J. Bang-Jensen and J. Huang, Kings in quasi-transitive digraphs, 
Discrete Math. 185 (1998) 19 - 27. 

\bibitem{msss} J. Bang-Jensen, J. Huang, and A. Yeo, Strongly connected spanning 
subdigraphs with the minimum number of arcs in quasi-transitive digraphs,
SIAM J. Discrete Math. 16 (2003) 335 - 343.  

\bibitem{galeana} H. Galeana-S\'anchez and R. Rojas-Monroy,
Kernels in quasi-transitive digraphs,
Discrete Math. 306 (2006) 1969 - 1974. 

\bibitem{golumbic} M.C. Golumbic, {\em Algorithmic Graph Theory and Perfect 
Graphs}, Academic Press, New York (1980).

\bibitem{hr} L. Haskins and D.J. Rose, Toward characterization of perfect 
elimination digraphs, SIAM J. Comput. 2 (1973) 217 - 224.

\bibitem{huang} J. Huang, On the structure of local tournaments,
J. Combin. Theory B (1995) 200 - 221.

\bibitem{meister} D. Meister and J.A. Telle, Chordal digraphs,
Theoret. Comput. Sci. 463 (2012) 73 - 83.

\bibitem{ye} Y.Y. Ye, On chordal digraphs and semi-strict chordal digraphs,
M.Sc. Thesis, University of Victoria, 2019. 

\end{thebibliography}
\end{document}